\documentclass{siamltex}
\pdfoutput=1
\usepackage{amsmath}
\usepackage{amsfonts}
\usepackage{graphicx}

\newcommand{\R}{\mathbb{R}}
\newcommand{\C}{\mathbb{C}}
\newcommand{\I}{\mathcal{I}}
\newcommand{\J}{\mathcal{J}}

\newcommand{\psd}{\succeq}
\newcommand{\pd}{\succ}
\newcommand{\nsd}{\preceq}

\DeclareMathOperator*{\tr}{\text{tr}}
\DeclareMathOperator*{\spn}{\text{span}}

\newcommand{\Ex}{\mathbb{E}}

\newcommand{\colsp}{\mathcal{R}}
\newcommand{\nullsp}{\mathcal{N}}
\newcommand{\elliptope}{\mathcal{E}}
\newcommand{\ones}{\mathbf{1}}
\newcommand{\face}{\mathcal{F}}

\newcommand{\blkdiag}{\text{blkdiag}}
\newcommand{\subspace}[1]{\mathcal{#1}}
\newcommand{\Sym}{\mathcal{S}}
\newcommand{\Amap}{\mathcal{A}}
\DeclareMathOperator*{\minimize}{\text{minimize}}
\DeclareMathOperator*{\maximize}{\text{maximize}}

\usepackage{subfig}
\usepackage{caption}

    \title{Diagonal and Low-Rank Matrix Decompositions,
            Correlation Matrices, and Ellipsoid Fitting\thanks{This research 
                was funded in part by Shell International Exploration and 
                Production, Inc.~under P.O.~450004440, and in part by the 
                Air Force Office of Scientific Research under 
                grant \#FA9550-11-1-0305. A preliminary version of parts of this 
                work appeared in the Master's thesis of the first-named author \cite{saunderson2011subspace}.}}

        \author{J.~Saunderson\footnotemark[2]\and V.~Chandrasekaran\footnotemark[2]\and P.~A.~Parrilo\footnotemark[2]\and
            A.~S.~Willsky\footnotemark[2]}

\begin{document}
    \maketitle
    \renewcommand{\thefootnote}{\fnsymbol{footnote}}
    \footnotetext[2]{Laboratory for Information and Decision Systems, Department of 
        Electrical Engineering and Computer Science, Massachusetts Institute of Technology, 
        Cambridge, MA 02139 (jamess@mit.edu, venkatc@mit.edu, parrilo@mit.edu, willsky@mit.edu).}

    \begin{abstract}
        In this paper we establish links between, and new results for, three problems that are
        not usually considered together. The first is a matrix decomposition problem that arises in
        areas such as statistical modeling and signal processing: given a matrix $X$ formed as the 
        sum of an unknown diagonal matrix and an unknown low rank positive semidefinite matrix, 
        decompose $X$ into these constituents.
        The second problem we consider is to determine the facial structure of the set of correlation 
        matrices, a convex set also known as the elliptope. This convex body, and particularly its facial structure,
        plays a role in applications from combinatorial optimization to mathematical finance.
        The third problem is a basic geometric question: given points $v_1,v_2,\ldots,v_n\in \R^k$
        (where $n > k$) determine whether there is a centered ellipsoid passing \emph{exactly} 
        through all of the points.

        We show that in a precise sense these three problems are equivalent. Furthermore 
        we establish a simple sufficient condition on a subspace $\subspace{U}$ 
        that ensures any positive semidefinite matrix $L$ with column space $\subspace{U}$
        can be recovered from $D+L$ for any diagonal matrix $D$ using a convex optimization-based
        heuristic known as minimum trace factor analysis. This result leads to a new understanding of
        the structure of rank-deficient correlation matrices and a simple condition on a set of 
        points that ensures there is a centered ellipsoid passing through them.
        \end{abstract}



        \pagestyle{myheadings}
        \thispagestyle{plain}
        \markboth{Diagonal and Low-Rank Matrix Decompositions}{Saunderson, Chandrasekaran, Parrilo, Willsky} 

\section{Introduction}
\label{sec:intro}

Decomposing a matrix as a sum of matrices with simple structure is a
fundamental operation with numerous applications. A matrix decomposition may
provide computational benefits, such as allowing the efficient solution of the
associated linear system in the square case. Furthermore, if the matrix arises from measurements
of a physical process (such as a sample covariance matrix), decomposing that
matrix can provide valuable insight about the structure of the physical
process.  

Among the most basic and well-studied additive matrix decompositions is the
decomposition of a matrix as the sum of a diagonal matrix and a low-rank
matrix. This decomposition problem arises in the factor analysis model in
statistics, which has been studied extensively since Spearman's original work
of 1904 \cite{spearman1904general}.  The same decomposition problem is known as
the Frisch scheme in the system identification literature
\cite{kalman1985identification}. For concreteness, in Section~\ref{sec:DOA} we
briefly discuss a stylized version of a problem in signal processing that under various assumptions
can be modeled as a (block) diagonal and low-rank decomposition problem.

Much of the literature on diagonal and low-rank matrix decompositions is in one
of two veins.  An early approach \cite{albert1944matrices} that has seen recent
renewed interest \cite{drton2007algebraic} is an algebraic one, where the
principal aim is to give a characterization of the vanishing ideal of the set of symmetric
$n\times n$ matrices that decompose as the sum of a diagonal matrix and a rank
$k$ matrix.  Such a characterization has only been obtained for the border
cases $k=1$, $k=n-1$ (due to Kalman \cite{kalman1985identification}), and the
recently resolved $k=2$ case (due to Brouwer and
Draisma~\cite{brouwer2011equivariant} following a conjecture by Drton et
al.~\cite{drton2007algebraic}). This approach does not (yet) offer scalable
algorithms for performing decompositions, rendering it unsuitable for many
applications including those in high-dimensional statistics, optics
\cite{fazel1998approximations}, and signal processing
\cite{saunderson2011subspace}. The other main approach to factor analysis is
via heuristic local optimization techniques, often based on the expectation
maximization (EM) algorithm \cite{dempster1977maximum}. This approach, while
computationally tractable, typically offers no provable performance guarantees.

A third way is offered by convex optimization-based methods for diagonal and
low-rank decompositions such as \emph{minimum trace factor
analysis} (MTFA), the idea and initial analysis of which dates at least to Ledermann's 1940 work
\cite{ledermann1940problem}. 
MTFA is computationally tractable, being based on a semidefinite program (see
Section~\ref{sec:bg-ps}), and yet offers the possibility of provable
performance guarantees. In this paper we provide a new analysis of MTFA that is
particularly suitable for high-dimensional problems. 

Semidefinite programming duality theory provides a link between this matrix
decomposition heuristic and the facial structure of the set of \emph{correlation matrices}--- 
positive semidefinite matrices with unit diagonal---also known as the 
\emph{elliptope}~\cite{laurent1995positive}.  This set is one of the simplest
of spectrahedra---affine sections of the positive semidefinite cone.
Spectrahedra are of particular interest for two reasons. First, spectrahedra are 
a rich class of convex sets that have many nice properties
(such as being facially exposed). Second, there are well-developed algorithms, efficient
both in theory and in practice, for optimizing linear functionals over spectrahedra. 
These optimization problems are known as semidefinite programs~\cite{vandenberghe1996semidefinite}.

The elliptope arises in semidefinite programming-based relaxations of
problems in areas such as combinatorial optimization (e.g.~the \textsc{max-cut}
problem~\cite{goemans1995improved}) and statistical mechanics (e.g.~the
$k$-vector spin glass problem~\cite{briet2010grothendieck}).   In addition, the
problem of projecting onto the set of (possibly low-rank) correlation matrices
has enjoyed considerable interest in mathematical finance and numerical
analysis in recent years~\cite{higham2002computing}. In each of these applications the
structure of the set of low-rank correlation matrices, i.e.~the facial
structure of this convex body, plays an important role.

Understanding the faces of the elliptope turns out to be related to the
following \emph{ellipsoid fitting problem}: given $n$ points in $\R^k$ (with $n > k$),
under what conditions on the points is there an ellipsoid centered at the
origin that passes \emph{exactly} through these points? While there is 
considerable literature on many ellipsoid-related problems, we are not aware
of any previous systematic investigation of this particular problem.

\subsection{Illustrative application: direction of arrival estimation}
\label{sec:DOA}
Direction of arrival estimation is a classical problem in
signal processing where (block) diagonal and low-rank decomposition
problems arise naturally. In this section we briefly discuss 
some stylized models of the direction of arrival estimation problem 
that can be reduced to matrix decomposition problems 
of the type considered in this paper. 

\begin{figure}
    \begin{center}
        \includegraphics[scale=1]{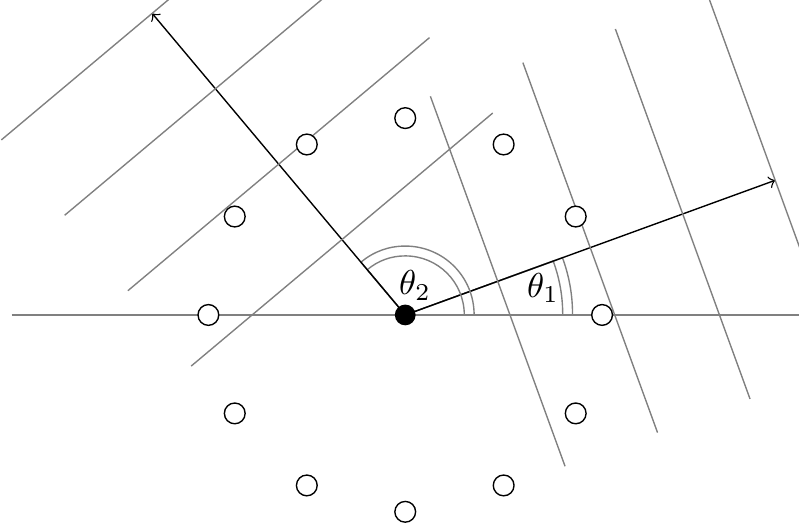}
    \end{center}
    \caption{\label{fig:doa} Plane waves from directions $\theta_1$ and $\theta_2$ arriving at an
        array of sensors equally spaced on a circle (a uniform circular array).}
\end{figure}

Suppose we have $n$ sensors at locations $(x_1,y_1),(x_2,y_2),\ldots,(x_n,y_n)
\in \R^2$ that are passively `listening' for waves (electromagnetic or
acoustic) at a known frequency from $r \ll n$ sources in the far field (so
that the waves are approximately plane waves when they reach the sensors). The
aim is to estimate the number of sources $r$ and their directions of arrival
$\theta = (\theta_1,\theta_2,\ldots,\theta_r)$ given sensor measurements and
knowledge of the sensor locations (see Figure~\ref{fig:doa}). 

A standard mathematical model for this problem (see \cite{krim1996two} for
a derivation) is to model the vector of sensor measurements $z(t)\in \C^n$ 
at time $t$ as 
\begin{equation}
\label{eq:array}
z(t) = A(\theta)s(t) + n(t)
\end{equation}
where $s(t)\in \C^r$ is the vector of baseband signal waveforms from the sources, 
$n(t)\in \C^n$ is the vector of sensor measurement noise, and 
$A(\theta)$ is the $n\times r$ matrix with complex entries 
$[A(\theta)]_{ij} = e^{-k\sqrt{-1}(x_i\cos(\theta_j)+y_i\sin(\theta_j))}$, 
with $k$ a positive constant related to the frequency of the waves being sensed. 

The column space of $A(\theta)$ contains all the information about the
directions of arrival $\theta$.  As such, subspace-based approaches to
direction of arrival estimation aim to estimate the column space of $A(\theta)$
(from which a number of standard techniques can be employed to estimate
$\theta$).

Typically $s(t)$ and $n(t)$ are modeled as zero-mean stationary white Gaussian
processes with covariances $\Ex[s(t)s(t)^H] = P$ and $\Ex[n(t)n(t)^H] = Q$
respectively (where $A^H$ denotes the Hermitian transpose of $A$ and
$\Ex[\cdot]$ the expectation). In the simplest setting, $s(t)$ and $n(t)$ are 
assumed to be uncorrelated so that the covariance of the sensor measurements at
any time is 
\[ \Sigma = A(\theta)PA(\theta)^H + Q.\]
The first term is Hermitian positive semidefinite with rank $r$, i.e.~the
number of sources. Under the assumption that spatially well-separated sensors (such as in a sensor network)
have uncorrelated measurement noise $Q$ is diagonal. In this case the covariance
$\Sigma$ of the sensor measurements decomposes as a sum of a positive
semidefinite matrix of rank $r\ll n$ and a diagonal matrix. Given an
approximation of $\Sigma$ (e.g.~a sample covariance) approximately performing
this diagonal and low-rank matrix decomposition allows the estimation of the
column space of $A(\theta)$ and in turn the directions of arrival.

A variation on this problem occurs if there are multiple sensors at each
location, sensing, for example, waves at different frequencies. Again under the assumption that 
well-separated sensors have uncorrelated measurement noise, and sensors at the 
same location have correlated measurement noise, the sensor noise covariance matrix
$Q$ would be \emph{block-diagonal}. As such the covariance of all of the
sensor measurements would decompose as the sum of a low-rank matrix (with rank
equal to the total number of sources over all measured frequencies) and
a block-diagonal matrix.

A block-diagonal and low-rank decomposition problem also arises if the second-order
statistics of the noise have certain \emph{symmetries}. This might occur in cases where
the sensors themselves are arranged in a symmetric way (such as in the uniform circular array shown in 
Figure~\ref{fig:doa}).
In this case there is a unitary matrix $T$ (depending only on the symmetry group of the array)
such that $TQT^H$ is \emph{block-diagonal} \cite{shah2011group}. Then the covariance of the sensor
measurements, when written in coordinates with respect to $T$, is 
\[T\Sigma T^H = TA(\theta)PA(\theta)^HT^H + TQT^H\]
which has a decomposition as the sum of a block diagonal matrix and a rank $r$ Hermitian positive semidefinite 
matrix (as conjugation by $T$ does not change the rank of this term).

Note that the matrix decomposition problems discussed in this section involve
Hermitian matrices with complex entries, rather than the symmetric matrices with real entries 
considered elsewhere in this paper. It is straightforward to generalize the main problems and results 
throughout the paper to the complex setting.

\subsection{Contributions}
\label{sec:contrib}
\paragraph{Relating MTFA, correlation matrices, and ellipsoid fitting}
We introduce and make explicit the links between the analysis of MTFA, the facial structure of
the elliptope, and the ellipsoid fitting problem, showing
that these problems are, in a precise sense, equivalent (see
Proposition~\ref{prop:links}).  As such, we relate a basic
problem in statistical modeling (tractable diagonal and low-rank matrix
decompositions), a basic problem in convex algebraic geometry (understanding
the facial structure of perhaps the simplest of spectrahedra), and a basic
geometric problem. 

\paragraph{A sufficient condition for the three problems}
The main result of the paper is to establish a new, simple, sufficient
condition on a subspace $\subspace{U}$ of $\R^n$ that ensures that MTFA correctly
decomposes matrices of the form $D^\star + L^\star$ where $\subspace{U}$ is the column
space of $L^\star$.  The condition is stated in terms of a measure
of~\emph{coherence} of a subspace (made precise in Definition~\ref{def:coh}).
Informally, the coherence of a subspace is a real
number between zero and one that measures how close the subspace is to
containing any of the elementary unit vectors. 
This result can be translated into new results for the other two problems
under consideration based on the relationship between the analysis of 
MTFA, the faces of the elliptope, and ellipsoid fitting. 

\paragraph{Block-diagonal and low-rank decompositions}

In Section~\ref{sec:block} we turn our attention to the \emph{block}-diagonal and
low-rank decomposition problem, showing how our results generalize to
that setting. Our arguments combine our results for the diagonal and low-rank decomposition
case with an understanding of the symmetries of the block-diagonal 
and low-rank decomposition problem. 

\subsection{Outline}
\label{sec:outline}

The remainder of the paper is organized as follows. 
We describe notation, give some background on semidefinite programming, and 
provide precise problem statements in Section~\ref{sec:bg-ps}.
In Section~\ref{sec:links} we present our first contribution by establishing
relationships between the success of MTFA, the faces of the elliptope, 
and ellipsoid fitting. We then illustrate these connections by noting the equivalence
of a known result about the faces of the elliptope, and a known result about MTFA, 
and translating these into the context of ellipsoid fitting.
Section~\ref{sec:suff} is focused on establishing and interpreting our main result: 
a sufficient condition for the three problems based on a coherence inequality. Finally in
Section~\ref{sec:block} we generalize our results to the analogous tractable
block-diagonal and low-rank decomposition problem.

\section{Background and problem statements}
\label{sec:bg-ps}

\subsection{Notation}
\label{sec:notation}
If $x,y\in \R^n$ we denote by $\langle x,y\rangle = \sum_{i=1}^{n} x_iy_i$ the
standard Euclidean inner product and by $\|x\|_2 = \langle x,x\rangle^{1/2}$
the corresponding Euclidean norm.  We write $x \geq 0$ and $x > 0$ to indicate
that $x$ is entry-wise non-negative and strictly positive, respectively.
Correspondingly, if $X,Y\in \Sym^n$, the set of $n \times n$ symmetric matrices, then we denote by
$\langle X,Y\rangle = \tr(XY)$ the trace inner product and by $\|X\|_F =
\langle X,X\rangle^{1/2}$ the Frobenius norm. We write $X \psd 0$ and $X \pd 0$ to indicate 
that $X$ is positive semidefinite and strictly positive definite, respectively. We write 
$\Sym_+^n$ for the cone of $n\times n$ positive semidefinite matrices.

The column space of a matrix $X$ 
is denoted $\colsp(X)$ and the nullspace is denoted $\nullsp(X)$. 
If $X$ is an $n\times n$ matrix then $\diag(X) \in \R^n$ is the diagonal of $X$.
If $x\in \R^n$ then $\diag^*(x)\in \Sym^n$ is the diagonal matrix with 
$[\diag^*(x)]_{ii} =x_i$ for $i=1,2,\ldots,n$. If $\subspace{U}$ is a subspace of $\R^n$ then
$P_{\subspace{U}}:\R^n\rightarrow \R^n$ denotes the orthogonal projector onto
$\subspace{U}$, that is the self-adjoint linear map such that
$\colsp(P_{\subspace{U}}) = \subspace{U}$, $P_{\subspace{U}}^2 = P_{\subspace{U}}$ and $\tr(P_{\subspace{U}}) =
\dim(\subspace{U})$.  

We use the notation $e_i$ for the vector with a one in the $i$th position and zeros elsewhere
and the notation $\ones$ to denote the vector all entries of which are one. We use the shorthand
$[n]$ for the set $\{1,2,\ldots,n\}$. The set of $n\times n$ correlation matrices, i.e.~positive semidefinite
matrices with unit diagonal, is denoted $\elliptope_n$. For brevity we typically refer to $\elliptope_n$ 
as the elliptope, and the elements of $\elliptope_n$ as correlation matrices.

\subsection{Semidefinite programming}
\label{sec:sdp}
The term semidefinite programming \cite{vandenberghe1996semidefinite} refers to 
convex optimization problems of the form
\begin{equation}
    \minimize_X\; \langle C,X \rangle \quad\text{subject to}\quad \left\{\begin{array}{rcl}  
        \Amap(X)\!\!\!\!& = &\!\!\!\! b\\
        X\!\!\!\! & \psd &\!\!\!\! 0\end{array}\right.\label{eq:sdp-primal}
\end{equation}
where $X$ and $C$ are $ n\times n$ symmetric matrices, $b\in \R^m$, and 
$\Amap : \Sym^{n} \rightarrow \R^m$ is a linear map. The dual semidefinite program 
is 
\begin{equation}
    \maximize_{y,S}\; \langle b,y\rangle \quad\text{subject to}\quad
    \left\{\begin{array}{rcl}
            C - \Amap^*(y) \!\!\!\! & = & \!\!\!\! S\\
            S \!\!\!\! & \psd & \!\!\!\! 0
        \end{array}\right.\label{eq:sdp-dual}
\end{equation}
where $\Amap^*:\R^m\rightarrow \Sym^n$ is the adjoint of $\Amap$. 

General semidefinite programs can be solved in polynomial time using interior 
point methods \cite{vandenberghe1996semidefinite}. While our focus in this paper is not on algorithms, 
we remark that for the structured semidefinite programs discussed in this paper, 
many different special-purpose methods have been devised.

The main result about semidefinite programming that we use is the following 
optimality condition (see~\cite{vandenberghe1996semidefinite} for example). 
\begin{theorem}
    \label{thm:sdp-opt}
    Suppose \eqref{eq:sdp-primal} and \eqref{eq:sdp-dual} are strictly feasible. Then 
    $X^\star$ and $(y^\star,S^\star)$ are optimal for the primal \eqref{eq:sdp-primal} and dual \eqref{eq:sdp-dual} respectively 
    if and only if $X^\star$ is primal feasible, $(y^\star,S^\star)$ is dual feasible and 
    $X^\star S^\star = 0$.
\end{theorem}

\subsection{Tractable diagonal and low-rank matrix decompositions}
\label{sec:dlr}

To decompose $X$ into a diagonal part and a positive semidefinite low-rank
part, we may try to solve the following rank minimization problem
\begin{equation*}
 	\minimize_{D,L} \;\rank(L)\quad\text{subject to}\quad
    \left\{\begin{array}{rcl}
            X \!\!\!\! & = & \!\!\!\! D + L\\
            L \!\!\!\! & \psd & \!\!\!\! 0\\
            D \!\!\!\! & & \!\!\!\!\!\!\!\!\text{diagonal.}\end{array}\right.
\end{equation*}
Since the rank function is non-convex and non-differentiable, it is not clear
how to solve this optimization problem directly. One approach that has been
successful for other rank minimization problems 
(for example those in \cite{mesbahi1997rank,recht2010guaranteed}), is to replace the rank function
with the trace function in the objective. This can be viewed as a convexification of the problem as
the trace function is the convex envelope of the rank function when restricted to positive semidefinite matrices
with spectral norm at most one. Performing this convexification leads to the semidefinite program we refer to as 
\emph{minimum trace factor analysis} (MTFA):
\begin{equation}
 	\minimize_{D,L} \;\tr(L)\quad\text{subject to}\quad
    \left\{\begin{array}{rcl}
            X \!\!\!\! & = & \!\!\!\! D + L\\
            L \!\!\!\! & \psd & \!\!\!\! 0\\
            D \!\!\!\! & & \!\!\!\!\!\!\!\!\text{diagonal.}\end{array}\right.\label{eq:mtfa0}
\end{equation}
It has been shown by Della Riccia and Shapiro \cite{della1982minimum} that if MTFA is feasible it has a unique
optimal solution. One central concern of this paper is
to understand when the diagonal and low-rank decomposition of a matrix given by MTFA is `correct' in the following sense.
\paragraph{Recovery problem I}
    Suppose $X$ is a matrix of the form $X = D^\star + L^\star$ where $D^\star$ is
diagonal and $L^\star$ is positive semidefinite. What conditions on 
$(D^\star,L^\star)$ ensure that $(D^\star,L^\star)$ is the unique 
optimum of MTFA with input $X$? 

We establish in Section~\ref{sec:links} that whether $(D^\star,L^\star)$ is the unique optimum of 
MTFA with input $X = D^\star + L^\star$ depends only on the column space of $L^\star$, motivating the
following definition.
\begin{definition}
    \label{def:rMTFA}
    A subspace $\subspace{U}$ of $\R^n$ is \emph{recoverable by MTFA} if for every diagonal $D^\star$ and 
    every positive semidefinite $L^\star$ with column space $\subspace{U}$, $(D^\star,L^\star)$ is the unique
    optimum of MTFA with input $X = D^\star + L^\star$. 
\end{definition}

In these terms, we can restate the recovery problem succinctly as follows.
\paragraph{Recovery problem II}
    Determine which subspaces of $\R^n$ are recoverable by MTFA.

Much of the basic analysis of MTFA, including optimality conditions and relations between 
minimum rank and minimum trace factor analysis, was carried out in a sequence of papers
by Shapiro \cite{shapiro1982rank,shapiro1982weighted,shapiro1985identifiability} and Della 
Riccia and Shapiro \cite{della1982minimum}. 
More recently, Chandrasekaran 
et al.~\cite{chandrasekaran2011rank} and Cand{\`e}s et al.~\cite{candes2011robust}
considered convex optimization heuristics for decomposing a matrix as a 
sum of a sparse and low-rank matrix. Since a diagonal matrix is certainly sparse, 
the analysis in \cite{chandrasekaran2011rank} can be specialized to give fairly conservative 
sufficient conditions for the success of MTFA.  

The diagonal and low-rank decomposition problem
can also be interpreted as a low-rank matrix completion problem, 
where we are given all the entries of a low-rank matrix except the diagonal, and 
aim to correctly reconstruct the diagonal entries. As such, this paper is closely 
related to the ideas and techniques used in the work of Cand{\`e}s and Recht \cite{candes2009exact}
and a number of subsequent papers on this topic. We would like to emphasize a key point of distinction
between that line of work and the present paper.
The recent low-rank matrix completion literature largely focuses on determining the 
proportion of \emph{randomly selected} entries of a low-rank matrix that need to be revealed to be able to reconstruct that low-rank matrix 
using a tractable algorithm. The results of this paper, on the other hand, can be interpreted as
attempting to understand which low-rank matrices can be reconstructed from a \emph{fixed} and 
quite canonical pattern of revealed entries.

\subsection{Faces of the elliptope}
\label{sec:elliptope}

The faces of the cone of $n\times n$ positive semidefinite matrices are all of
the form
\begin{equation} 
    \label{eq:psd-faces}
    \face_{\subspace{U}} = \{X \psd 0: \nullsp(X) \supseteq \subspace{U}\}
\end{equation}
where $\subspace{U}$ is a subspace of $\R^n$ \cite{laurent1995positive}. Conversely given any subspace 
$\subspace{U}$ of $\R^n$, $\face_{\subspace{U}}$ is a face of $\Sym_+^n$.  As a
consequence, the faces of $\elliptope_n$ are all of the form
\begin{equation}
    \label{eq:elliptope}
    \elliptope_n \cap \face_{\subspace{U}} = \{X \psd 0: \nullsp(X) \supseteq \subspace{U},\; \diag(X) = \ones\}
\end{equation}
where $\subspace{U}$ is a subspace of $\R^n$ \cite{laurent1995positive}. It is \emph{not} the case, however, that for
every subspace $\subspace{U}$ of $\R^n$ there is a correlation matrix with nullspace
containing $\subspace{U}$, motivating the following definition.
\begin{definition}[\cite{laurent1995positive}]
    \label{def:realizable}
    A subspace $\subspace{U}$ of $\R^n$ is \emph{realizable} if there is an $n\times n$
    correlation matrix $Q$ such that $\nullsp(Q) \supseteq \subspace{U}$.
\end{definition}

The problem of understanding the facial structure of the set of correlation matrices 
can be restated as follows.
\paragraph{Facial structure problem} Determine which subspaces of $\R^n$ are realizable.

Much is already known about the faces of the elliptope.
For example, all possible dimensions of faces as well as 
polyhedral faces, are known \cite{laurent1995facial}. Characterizations of the realizable 
subspaces of $\R^n$ of dimension $1$, $n-2$, and $n-1$ are given in \cite{delorme1993combinatorial} and 
implicitly in \cite{laurent1995positive} and \cite{laurent1995facial}. Nevertheless, little 
is known about which $k$ dimensional subspaces of $\R^n$ are realizable for general $n$ and $k$.

\subsection{Ellipsoid fitting} 
\label{sec:ef} 
\paragraph{Ellipsoid fitting problem I} What conditions on a collection of $n$ points in $\R^k$
ensure that there is a centered ellipsoid passing \emph{exactly} through all those points?

Let us consider some basic properties of this problem.
\paragraph{Number of points} If $n\leq k$ we can always fit an ellipsoid to the points. Indeed if
   $V$ is the matrix with columns $v_1,v_2,\ldots,v_n$ then the image of the unit sphere in 
   $\R^n$ under $V$ is a centered ellipsoid passing through $v_1,v_2,\ldots,v_n$. 
   If $n > \binom{k+1}{2}$ and the points are `generic' then we cannot fit a centered ellipsoid to them.
  This is because if we represent the ellipsoid by a symmetric $k\times k$ matrix $M$, the condition
  that it passes through the points (ignoring the positivity condition on $M$) means
  that $M$ must satisfy $n$ linearly independent equations.

\paragraph{Invariances} If $T\in GL(k)$ is an invertible linear map then there is 
an ellipsoid passing through $v_1,v_2,\ldots,v_n$ if and only if there is an ellipsoid
passing through $Tv_1,Tv_2,\ldots,Tv_n$. This means that whether there is an ellipsoid passing through
$n$ points in $\R^k$ does not depend on the actual set of $n$ points, but on a subspace of 
$\R^n$ related to the points. We summarize this observation in the following lemma.
\begin{lemma}
    \label{lem:efp}
    Suppose $V$ is a $k\times n$ matrix with row space $\subspace{V}$. If there is a centered ellipsoid in $\R^k$ 
    passing through the columns of $V$ then there is a centered ellipsoid passing through the columns of 
    any matrix $\tilde{V}$ with row space $\subspace{V}$.
\end{lemma}

Lemma~\ref{lem:efp} asserts that whether it is possible to fit an ellipsoid to $v_1,v_2,\ldots,v_n$ depends only on the 
row space of the matrix with columns given by the $v_i$, motivating the following definition. 
\begin{definition}
    \label{def:efp}
    A subspace $\subspace{V}$ of $\R^n$ has the \emph{ellipsoid fitting property} if there is 
    a $k\times n$ matrix $V$ with row space $\subspace{V}$ and a centered ellipsoid in $\R^k$
    that passes through each column of $V$.
\end{definition}

As such we can restate the ellipsoid fitting problem as follows.
\paragraph{Ellipsoid fitting problem II} Determine which subspaces of $\R^n$
have the ellipsoid fitting property.

\section{Relating ellipsoid fitting, diagonal and low-rank decompositions, and correlation matrices}
\label{sec:links}

In this section we show that the ellipsoid fitting problem, the recovery problem, and the facial structure
problem are equivalent in the following sense.
\begin{proposition}
    \label{prop:links}
    Let $\subspace{U}$ be a subspace of $\R^n$.
    Then the following are equivalent:
    \begin{remunerate}
        \item \label{item:rMTFA} $\subspace{U}$ is recoverable by MTFA.
        \item \label{item:realizable} $\subspace{U}$ is realizable.
        \item \label{item:efp} $\subspace{U}^\perp$ has the ellipsoid fitting property.
    \end{remunerate}
\end{proposition}
\begin{proof}
    To see that~\ref{item:realizable} implies~\ref{item:efp}, 
    let $V$ be a $k\times n$ matrix with nullspace $\subspace{U}$ and let $v_i$ denote the $i$th 
column of $V$. If $\subspace{U}$ is realizable there is a correlation matrix $Y$ with nullspace 
containing $\subspace{U}$. Hence there is some $M\psd 0$ such that $Y = V^T M V$ and $v_i^TMv_i = 1$ for 
$i\in[n]$. Since $V$ has nullspace $\subspace{U}$, it has row space $\subspace{U}^\perp$. Hence the subspace
$\subspace{U}^\perp$ has the ellipsoid fitting property. By reversing the argument we see that the converse 
also holds.

The equivalence of~\ref{item:rMTFA} and~\ref{item:realizable} arises from semidefinite 
programming duality. Following a slight reformulation, MTFA \eqref{eq:mtfa0} 
can be expressed as
\begin{equation}
    \label{eq:mtfa-primal}
    \maximize_{d,L} \; \langle \ones,d\rangle \quad\text{subject to}\quad
    \left\{\begin{array}{rcl}
            X\!\!\!\! & = & \!\!\!\!\diag^*(d) + L\\
            L \!\!\!\! & \psd  & \!\!\!\!0\end{array}\right.
\end{equation}
and its dual as
\begin{equation}
    \label{eq:mtfa-dual}
    \minimize_{Y} \;\langle X,Y\rangle \quad\text{subject to}\quad 
    \left\{\begin{array}{rcl}
            \diag(Y)\!\!\!\! & = & \!\!\!\!\ones\\
            Y \!\!\!\! & \psd & \!\!\!\! 0\end{array}\right.
\end{equation}
which is clearly just the optimization of the linear functional defined by $X$ over the elliptope.
We note that~\eqref{eq:mtfa-primal} is exactly in the standard dual form~\eqref{eq:sdp-dual}
for semidefinite programming and correspondingly that~\eqref{eq:mtfa-dual} is in the standard
primal form~\eqref{eq:sdp-primal} for semidefinite programming.

Suppose $\subspace{U}$ is recoverable by MTFA\@. Fix a diagonal matrix $D^\star$ and 
a positive semidefinite matrix $L^\star$ with column space $\subspace{U}$ and let
$X = D^\star + L^\star$. Since~\eqref{eq:mtfa-primal} 
and \eqref{eq:mtfa-dual} are strictly feasible, by 
Theorem~\ref{thm:sdp-opt} (optimality conditions for semidefinite programming),
the pair $(\diag(D^\star),L^\star)$ is an optimum of~\eqref{eq:mtfa-primal} if and only
if there is some correlation matrix $Y^\star$ such that $Y^\star L^\star = 0$.
Since $\colsp(L^\star) = \subspace{U}$ this implies that $\subspace{U}$ is realizable. 
Conversely, if $\subspace{U}$ is realizable, there is some $Y^\star$ such that 
$Y^\star L^\star = 0$ for every $L^\star$ with column space $\subspace{U}$, showing that 
$\subspace{U}$ is recoverable by MTFA.
\end{proof}

\paragraph{Remark} We note that in the proof of Proposition~\ref{prop:links} we established
that the two versions of the recovery problem stated in Section~\ref{sec:dlr} are actually 
equivalent. In particular, whether $(D^\star,L^\star)$ is the optimum of MTFA with input
$X = D^\star + L^\star$ depends only on the column space of $L^\star$.

\subsection{Certificates of failure}
\label{sec:failure}
We can prove that a subspace $\subspace{U}$ is realizable by constructing a correlation matrix
with nullspace containing $\subspace{U}$. We can prove that a subspace is \emph{not} realizable 
by constructing a matrix that \emph{certifies}
this fact. Geometrically, a subspace $\subspace{U}$ is realizable 
if and only if the subspace $\subspace{L}_{\subspace{U}} = \{X\in \Sym^n: \nullsp(X) \supseteq \subspace{U}\}$ of symmetric matrices 
intersects with the elliptope. So a certificate that $\subspace{U}$ is not realizable is a hyperplane
in the space of symmetric matrices that strictly separates the elliptope from $\subspace{L}_{\subspace{U}}$. 
The following lemma describes the structure of these separating hyperplanes.
\begin{lemma}
    \label{lem:failure}
    A subspace $\subspace{U}$ of $\R^n$ is \emph{not} realizable if and only if 
    there is a diagonal matrix $D$ such that $\tr(D) > 0$ and 
    $v^TDv \leq 0$ for all $v\in \subspace{U}^\perp$. 
\end{lemma}

\begin{proof}
    By Proposition~\ref{prop:links}, $\subspace{U}$ is not realizable if and only if 
    $\subspace{U}^\perp$ does not have the ellipsoid fitting property. Let $V$ be a 
    $k \times n$ matrix with row space $\subspace{U}^\perp$. Then $\subspace{U}^\perp$
    does not have the ellipsoid fitting property if and only if we cannot find an 
    ellipsoid passing through the columns of $V$, i.e.~the semidefinite
    program 
    \begin{equation}
        \label{eq:ef-feas}
        \minimize_{M} \; \langle 0,M\rangle \quad\text{subject to}\quad
        \left\{\begin{array}{rcl} \diag(V^TMV)\!\!\!\! & = & \!\!\!\!\ones\\
                M \!\!\!\! & \psd & \!\!\!\! 0 \end{array}\right.
        \end{equation}
        is infeasible. The semidefinite programming dual of \eqref{eq:ef-feas}
        is 
    \begin{equation}
        \label{eq:ef-feas-dual}
        \maximize_d \; \langle d,\ones\rangle \quad\text{subject to}\quad
        \left\{\begin{array}{rcl} V\diag^*(d)V^T\!\!\!\! & \nsd & \!\!\!\! 0.\end{array}\right.
        \end{equation}
        Since \eqref{eq:ef-feas-dual} is clearly always feasible, by strong duality (which holds
        because both primal and dual problems are strictly feasible)
        \eqref{eq:ef-feas} is infeasible if and only if \eqref{eq:ef-feas-dual} is unbounded. This occurs 
        if and only if there is some $d$ with $\sum_{i\in [n]} d_i > 0$ and yet 
        $V\diag^*(d)V^T \nsd 0$. Then $D = \diag^*(d)$ has the properties in the statement of the lemma.
\end{proof}

\subsection{Exploiting connections: results for one dimensional subspaces}
\label{sec:necc}

In 1940, Ledermann~\cite{ledermann1940problem} characterized the one dimensional 
subspaces that are recoverable by MTFA\@. In 1990, Grone et al.~\cite{grone1990extremal} gave a 
necessary condition for a subspace to be realizable. In 1993, independently of Ledermann's work, 
Delorme and Poljak~\cite{delorme1993combinatorial} showed that this condition is also 
sufficient for one dimensional subspaces. Since we have established 
that a subspace is recoverable by MTFA if and only if it is realizable, 
Ledermann's result and Delorme and Poljak's results are equivalent. In this section 
we translate these equivalent results into the context of the ellipsoid fitting problem, 
giving a geometric characterization of when it is possible to fit a centered ellipsoid 
to $k+1$ points in $\R^k$.

Delorme and Poljak state their result in terms of the following definition.
\begin{definition}[\cite{delorme1993combinatorial}]
    \label{def:bal}
    A vector $u\in \R^n$ is \emph{balanced} if, for all $i\in[n]$,
    \begin{equation}
        \label{eq:balance}
        |u_i| \leq \sum_{j\neq i} |u_j|.
    \end{equation}
    If the inequality is strict we say that $u$ is \emph{strictly balanced}.
\end{definition}

In the following, the necessary condition is due to Grone et al.~\cite{grone1990extremal} 
and the sufficient condition is due to Ledermann \cite{ledermann1940problem} (in the context of 
the analysis of MTFA) and Delorme and Poljak \cite{delorme1993combinatorial} (in the context of 
the facial structure of the elliptope). We state the result only in terms of 
realizability of a subspace.
\begin{theorem}
    \label{thm:DP}
    If a subspace $\subspace{U}$ of $\R^n$ is realizable 
    then every $u\in \subspace{U}$ is balanced. 
    If $\subspace{U} = \spn\{u\}$ is one-dimensional then $\subspace{U}$ is realizable 
    if and only if $u$ is balanced. 
\end{theorem}

The balance condition has a particularly natural geometric interpretation
in the ellipsoid fitting setting (Lemma~\ref{lem:bal-ef}, below). The proof is a fairly 
straightforward application of linear programming duality, which we defer to
Appendix~\ref{app:pfs}.
\begin{lemma}
    \label{lem:bal-ef}
    Suppose $V$ is any $k \times n$ matrix with $\nullsp(V) = \subspace{U}$. 
   Denote the columns of $V$ by $v_1,v_2,\ldots,v_n \in \R^{k}$. 
    Then every $u\in \subspace{U}$ is balanced if and only if for each 
    $i\in [n]$, $v_i$ lies on the boundary of the convex hull of 
    $\pm v_1,\pm v_2,\ldots,\pm v_n$.
\end{lemma}

By combining Theorem~\ref{thm:DP} with Lemma~\ref{lem:bal-ef}, 
we are in a position to interpret Theorem~\ref{thm:DP} 
purely in terms of ellipsoid fitting.
\begin{corollary}
    \label{cor:bal-ef}
    If there is an ellipsoid passing through $\pm v_1,\pm v_2,\ldots,\pm v_n\in \R^k$ then
    $\pm v_1,\pm v_2,\ldots,\pm v_n$ lie on the boundary of their convex hull.  If, in addition,
    $k=n-1$ the converse also holds.
\end{corollary}

We note that $\pm v_1,\pm v_2,\ldots,\pm v_n$ lie on the boundary of their convex hull
if and only if there exists \emph{some} convex set with boundary containing $\pm v_1,\pm v_2,\ldots,\pm v_n$.
In this geometric setting, it is clear that this is a necessary condition to be able to 
find a centered ellipsoid passing through the points, but not so obvious that it is 
sufficient if $k=n-1$. 

\section{A sufficient condition for the three problems}
\label{sec:suff}
In this section we establish a new sufficient condition for a subspace $\subspace{U}$ 
of $\R^n$ to be realizable and consequently a sufficient condition for $\subspace{U}$ to 
be recoverable by MTFA and $\subspace{U}^\perp$ to have the ellipsoid fitting property.
Our condition is based on a simple property of a subspace known as coherence.

Given a subspace $\subspace{U}$ of $\R^n$, the coherence of $\subspace{U}$ is a measure of 
    how close the subspace is to containing any of the elementary unit 
    vectors. This notion was introduced (with a different scaling) by Cand{\`e}s and Recht in 
    their work on low-rank matrix completion \cite{candes2009exact}, although related quantities
    have played an important role in the analysis of sparse reconstruction problems since
    the work of Donoho and Huo \cite{donoho2001uncertainty}. 
    \begin{definition}
    \label{def:coh}
        If $\subspace{U}$ is a subspace of $\R^n$ then the \emph{coherence} of $\subspace{U}$ is
        \[ \mu(\subspace{U}) = \max_{i\in [n]}\|P_\subspace{U}e_i\|_2^2.\]
       % = \max_{i\in [n]} [P_U]_{ii}.\] put this in the appendix
    \end{definition}
    
    A basic property of coherence is that it satisfies the inequality
    \begin{equation}
        \label{eq:coh-ineq}
        \frac{\dim(\subspace{U})}{n} \leq \mu(\subspace{U}) \leq 1
    \end{equation}
    for any subspace $\subspace{U}$ of $\R^n$ \cite{candes2009exact}. This inequality, together with the definition
    of coherence, provides useful intuition about the properties of subspaces with low coherence,
    that is \emph{in}coherence. Any subspace with low coherence is necessarily
    of low dimension and far from containing any of the elementary unit vectors $e_i$.
    As such, any symmetric matrix with incoherent row/column spaces is necessarily of
    low-rank and quite different from being a diagonal matrix.

    \subsection{Coherence-threshold-type sufficient conditions}
    In this section we focus on finding the largest possible $\alpha$ such that 
\[ \text{$\mu(\subspace{U}) < \alpha \implies \subspace{U}$ is realizable,}\]
that is finding the best possible coherence-threshold-type sufficient condition for a subspace to be realizable.
Such conditions are of particular interest because the dependence they have 
on the ambient dimension and the dimension of the subspace is only the mild dependence implied by \eqref{eq:coh-ineq}. 
In contrast, existing results (e.g.~\cite{delorme1993combinatorial,laurent1995facial,laurent1995positive})
about realizability of subspaces hold only for specific combinations of the 
ambient dimension and the dimension of the subspace.

The following theorem, our main result, gives a sufficient condition for realizability based on a coherence-threshold 
condition. Furthermore, it establishes that this is the best possible coherence-threshold-type sufficient condition.

\begin{theorem}
    \label{thm:main}
    If $\subspace{U}$ is a subspace of $\R^n$ and $\mu(\subspace{U}) < 1/2$ then $\subspace{U}$ is realizable.
    On the other hand, given any $\alpha > 1/2$, there is a subspace $\subspace{U}$ with $\mu(\subspace{U}) = \alpha$ that is not realizable.
\end{theorem}
\begin{proof}
    We give the main idea of the proof, deferring some details to Appendix~\ref{app:pfs}. 
    Instead of proving that there is some $Y\in \mathcal{F}_{\subspace{U}} = \{Y \psd 0: \nullsp(Y) \supseteq\subspace{U}\}$ 
    such that $Y_{ii} = 1$ for $i\in [n]$, it suffices to choose a convex cone $\mathcal{K}$ that is an inner 
    approximation to $\mathcal{F}_{\subspace{U}}$ and establish that there is some $Y\in \mathcal{K}$ such that 
    $Y_{ii} = 1$ for $i\in [n]$. 
    One natural choice is to take $\mathcal{K} = \{P_{\subspace{U}^\perp}\diag^*(\lambda)P_{\subspace{U}^\perp}: \lambda\geq 0\}$, 
    which is clearly contained in $\mathcal{F}_{\subspace{U}}$. Note that there is some
    $Y \in \mathcal{K}$ such that $Y_{ii} = 1$ for all $i\in [n]$ if
    and only if there is $\lambda \geq 0$ such that 
    \begin{equation}
        \label{eq:lin-sys}
        \diag\left(P_{\subspace{U}^\perp}\diag^*(\lambda)P_{\subspace{U}^\perp}\right) = \ones.
    \end{equation}
    The rest of the proof of the sufficient condition involves showing that if $\mu(\subspace{U}) < 1/2$ then such a 
    non-negative $\lambda$ exists. We establish this in Lemma~\ref{lem:linsys}.

    Now let us construct, for any $\alpha > 1/2$, a subspace with coherence $\alpha$ that is not realizable.
    Let $\subspace{U}$ to be the subspace of $\R^2$ spanned by $u = (\sqrt{\alpha},\sqrt{1-\alpha})$. Then 
    $\mu(\subspace{U}) = \max\{\alpha,1-\alpha\} = \alpha$ and yet by Theorem~\ref{thm:DP}, 
    $\mathcal{U}$ is not realizable because $u$ is not balanced.
\end{proof}

\paragraph{Remarks} Theorem~\ref{thm:main} illustrates both the power and limitations of coherence-threshold-type 
conditions. On the one hand, since coherence is quite a coarse property of a subspace, the result applies
to `many' subspaces (see Proposition~\ref{prop:rand} in Section~\ref{sec:examples}). On the other hand, since coherence has
very mild dimension dependence, the power of coherence-threshold-type conditions is limited to their 
specialization to low-dimensional situations, such as one dimensional subspaces of $\R^2$.

\subsection{Interpretations of Theorem~\ref{thm:main}}
\label{sec:interp}

We now establish two corollaries of our coherence-threshold-type sufficient condition for realizability. These 
corollaries can be thought of as re-interpretations of the coherence inequality $\mu(\subspace{U}) < 1/2$
in terms of other natural quantities.

\paragraph{An ellipsoid-fitting interpretation}

With the aid of Proposition~\ref{prop:links} we reinterpret our coherence-threshold-type sufficient condition as
a sufficient condition on a set of points in $\R^k$ that ensures there is a centered ellipsoid passing through them.
The condition involves `sandwiching' the points between two ellipsoids (that depend on the points). Indeed, 
given $v_1,v_2,\ldots,v_n\in \R^k$ and $0 < \beta < 1$ we define the ellipsoid 
\[ \mathcal{E}_{\beta}(v_1,\ldots,v_n) = \text{$\{x\in \R^k: x^T(\textstyle\sum_{j=1}^n v_jv_j^T)^{-1}x \leq \beta\}$}.\]
\begin{definition}
    Given $0< \beta < 1$ the points $v_1,v_2,\ldots,v_n$ satisfy the \emph{$\beta$-sandwich condition} if 
    \[\{v_1,v_2,\ldots,v_n\} \subset \mathcal{E}_1(v_1,\ldots,v_n) \setminus \mathcal{E}_\beta(v_1,\ldots,v_n).
        %= \{x\in \R^k:\alpha < x^T(\textstyle\sum_{j=1}^{n}v_jv_j^T)^{-1}x \leq 1.
    \]
\end{definition}

\begin{figure}
    \begin{center}
        \includegraphics{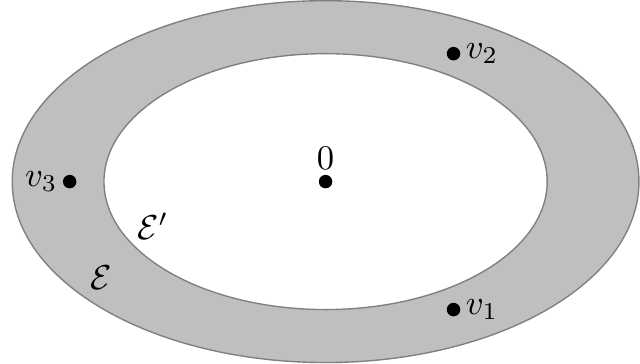}
    \end{center}
    \caption{\label{fig:sandwich} The ellipsoids 
        shown are $\mathcal{E} = \mathcal{E}_1(v_1,v_2,v_3)$ and $\mathcal{E}' = \mathcal{E}_{1/2}(v_1,v_2,v_3)$.
        There is an ellipsoid passing through $v_1,v_2$ and $v_3$ because the points are 
    sandwiched between $\mathcal{E}$ and $\mathcal{E}'$.}
\end{figure}

The intuition behind this definition (illustrated in Figure~\ref{fig:sandwich}) 
is that if the points satisfy the $\beta$-sandwich condition for $\beta$ close
to one, then they are confined to a thin elliptical shell that is adapted to their position. One might expect that 
it is `easier' to fit an ellipsoid to points that are confined in this way. Indeed this is the case.

\begin{corollary}
    \label{cor:ef}
    If $v_1,v_2,\ldots,v_n\in \R^k$ satisfy the $1/2$-sandwich condition then 
    then there is a centered ellipsoid passing through $v_1,v_2,\ldots,v_n$. 
\end{corollary}

\begin{proof}
    Let $V$ be the $k\times n$ matrix with columns given by the $v_i$, and let $\subspace{U}$ be the nullspace of $V$. 
    Then the orthogonal projection onto the row space of $V$ is $P_{\subspace{U}^\perp}$, and can be written as
    \[ P_{\subspace{U}^\perp} = V^T(VV^T)^{-1}V.\]
    Our assumption that the points satisfy the $1/2$-sandwich condition is 
    equivalent to assuming that $1/2 < [P_{\subspace{U}^\perp}]_{ii} \leq 1$ for all $i\in [n]$ 
    or alternatively that \[
        \mu(\subspace{U}) = \max_{i\in [n]} [P_{\subspace{U}}]_{ii} = 1-\min_{i\in [n]}[P_{\subspace{U}^\perp}]_{ii} < 1/2.
    \]
    From Theorem~\ref{thm:main} we know that $\mu(\subspace{U}) < 1/2$ implies that $\subspace{U}$ is realizable. Invoking 
    Proposition~\ref{prop:links} we then conclude that there is a centered ellipsoid passing through $v_1,v_2,\ldots,v_n$.    
\end{proof}

\paragraph{A balance interpretation}
In Section~\ref{sec:necc} we saw that if a subspace $\subspace{U}$ is
realizable, every $u\in \subspace{U}$ is balanced. The sufficient condition of 
Theorem~\ref{thm:main} can be expressed in terms of a balance condition on 
the element-wise square of the elements of a subspace. 
(In what follows $u\circ u$ denotes the element-wise square of a vector in $\R^n$.)
    \begin{corollary}
        \label{cor:bal-suff}
        Suppose $\subspace{U}$ is a subspace of $\R^n$. If $u\circ u$ is strictly balanced
        for every $u\in \subspace{U}$ then $\subspace{U}$ is realizable.
    \end{corollary}
    \begin{proof}
        It suffices to show that if for every $u\in \subspace{U}$, $u\circ u$ is strictly balanced,
        then $\mu(\subspace{U}) <1/2$ (although we could reverse the argument to establish the equivalence
        of these conditions). If $u\circ u$ is strictly balanced for all $u\in \subspace{U}$ then
        for all $i\in [n]$ and all $u\in \subspace{U}$
        \begin{equation}
            \label{eq:sqbal}
            2\langle e_i,u\rangle^2 < \sum_{j=1}^{n}\langle e_j,u\rangle^2 = \|u\|_2^2.
        \end{equation}
       Since $\|P_{\subspace{U}}e_i\|_2 = \max_{u\in \subspace{U}\setminus\{0\}} \langle e_i,u\rangle/\|u\|_2$,
       it follows from \eqref{eq:sqbal} that 
       $2\|P_{\subspace{U}}e_i\|_2^2 < 1$. Since this holds for all $i\in [n]$ it follows that $\mu(\subspace{U}) < 1/2$.
   \end{proof}
       
   \paragraph{Remark} Suppose $\subspace{U} = \spn\{u\}$ is a one-dimensional subspace of $\R^n$.
   We have just established that if $u\circ u$ is strictly balanced then $\subspace{U}$ is 
   realizable and so (by Theorem~\ref{thm:DP}) $u$ must be balanced. We note that it is straightforward
   to establish directly that if $u\circ u$ is balanced then $u$ is balanced by
   using the definition of balance and the fact that $\|x\|_1 \geq \|x\|_2$ for any $x\in \R^n$. 

\subsection{Examples}
\label{sec:examples}
To gain more intuition for what Theorem~\ref{thm:main} means, we consider its implications 
in two particular cases. First, we compare the characterization of when it is possible 
to fit an ellipsoid to $k+1$ points in $\R^k$ (Corollary~\ref{cor:bal-ef}) with the specialization of our sufficient 
condition to this case (Corollary~\ref{cor:ef}). This comparison provides some insight into how conservative our sufficient 
condition is. Second, we investigate the coherence properties of suitably `random' subspaces. This provides 
intuition about whether or not $\mu(\subspace{U}) < 1/2$ is a very restrictive condition. In particular, 
we establish that `most' subspaces of $\R^n$ with dimension bounded above by $(1/2 - \epsilon)n$ are realizable.

    \paragraph{Fitting an ellipsoid to $k+1$ points in $\R^k$}
    Recall that Ledermann and Delorme and Poljak's result, interpreted in terms of ellipsoid fitting,
    tells us that we can fit an ellipsoid to $k+1$ points
    $v_1,\ldots,v_{k+1}\in \R^k$ if and only if those points are on the
    boundary of the convex hull of $\{\pm v_1,\ldots,\pm v_{k+1}\}$ (see 
    Corollary~\ref{cor:bal-ef}). We now compare this characterization with 
    the $1/2$-sandwich condition, which is sufficient by Corollary~\ref{cor:ef}. 
    
\begin{figure}
        \begin{center}
            \subfloat[][The shaded set is $R$, those points $v$ for which we can fit an ellipsoid through 
            $v$ and the standard basis vectors.]{\includegraphics{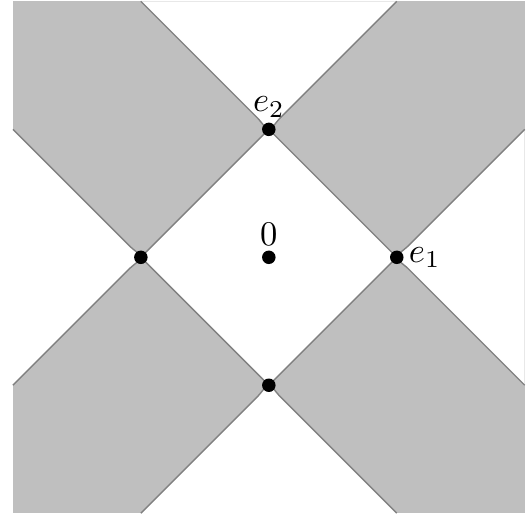}\label{fig:3ptsneccsuff}}\hfill
            \subfloat[][The shaded set is $R'$, those points $v$ such that $v,e_1$ and $e_2$ satisify the condition 
            of Corollary~\ref{cor:ef}.]{\includegraphics{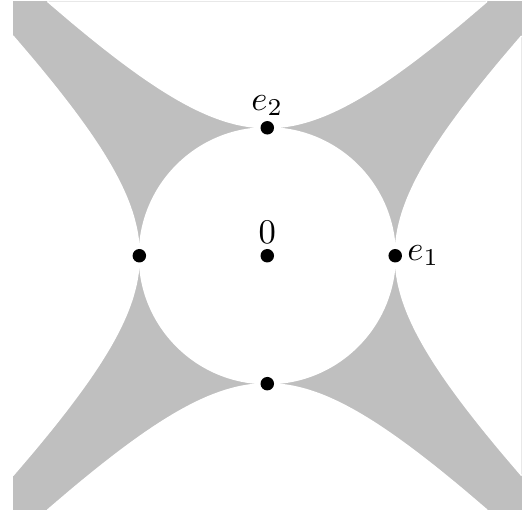}\label{fig:3ptssuff}}
        \end{center}
        \caption{\label{fig:3pts} Comparing our sufficient condition for ellipsoid fitting (Corollary~\ref{cor:ef}) with 
            the characterization (Corollary~\ref{cor:bal-ef}) in the case of fitting an ellipsoid to $k+1$ points in $\R^k$.}
    \end{figure}

    Without loss of generality we assume that $k$ of the
    points are $e_1,\ldots,e_k$, the standard basis vectors, and compare the
    conditions by considering the set of locations of the $k+1$st point $v\in \R^k$ 
    for which we can fit an ellipsoid through all $k+1$ points. Corollary~\ref{cor:bal-ef} 
    gives a characterization of this region as
    \[ R = \{v\in \R^k: \sum_{j=1}^k|v_j| \geq 1, \; |v_i| - \sum_{j\neq i}|v_j| \leq 1\quad\text{for $i\in[k]$}\}\]
    which is shown in Figure~\ref{fig:3ptsneccsuff} in the case $k=2$. The set
    of $v$ such that $v,e_1,\ldots,e_n$ satisfy the $1/2$-sandwich condition
    can be written as
    \begin{align*} 
        R' & = \{v\in \R^k: v^T(I + vv^T)^{-1}v > 1/2,\; e_i^T(I+vv^T)^{-1}e_i > 1/2 \quad\text{for $i\in[k]$}\}\\
        & = \{v\in \R^k: \sum_{j=1}^{k}v_j^2 > 1,\; v_i^2 - \sum_{j\neq i}v_j^2 < 1 \quad\text{for $i\in[k]$}\}
    \end{align*}
     which is shown in Figure~\ref{fig:3ptssuff}. It is clear that $R'
     \subseteq R$.

    \paragraph{Realizability of random subspaces}
    Suppose $\subspace{U}$ is a subspace generated by taking the column space of an
    $n\times r$ matrix with i.i.d.~standard Gaussian entries. For what values of $r$ 
    and $n$ does such a subspace have $\mu(\subspace{U}) < 1/2$ with high probability, 
    i.e.~satisfy our sufficient condition for being realizable? 
    
    The following result essentially shows that 
    for large $n$, `most' subspaces of dimension at most $(1/2-\epsilon)n$
    are realizable. This suggests that MTFA is a very good heuristic for diagonal and 
    low-rank decomposition problems in the high-dimensional setting. Indeed  
    `most' subspaces of dimension up to one half the 
    ambient dimension---hardly just low-dimensional subspaces---are recoverable by MTFA. 
    \begin{proposition}
        \label{prop:rand}
        Let $0 < \epsilon < 1/2$ be a constant and suppose $n > 6/(\epsilon^2-2\epsilon^3)$. 
        There are positive constants $\bar{c}$, $\tilde{c}$, (depending only on $\epsilon$) such that if 
        $\subspace{U}$ is a random $(1/2 - \epsilon)n$ dimensional subspace of $\R^n$ 
        then
        \[ \Pr[\text{$\subspace{U}$ is realizable}] \geq 1-\bar{c}\sqrt{n}e^{-\tilde{c} n}.\]
    \end{proposition}

    We provide a proof of this result in Appendix~\ref{app:pfs}. The main idea is that the coherence 
    of a random $r$ dimensional subspace of $\R^n$ is
    the maximum of $n$ random variables that concentrate around their mean of $r/n$ for large $n$.

    To illustrate the result, we consider the case where $\epsilon = 1/4$ and $n>192$. 
    Then (by examining the proof in Appendix~\ref{app:pfs}) we see that we can take
    $\tilde{c} = 1/24$ and $\bar{c} = 24/\sqrt{3\pi} \approx 7.8$. Hence if $n > 192$ and $\subspace{U}$ is 
    a random $n/4$ dimensional subspace of $\R^n$ we have that
    \[ \Pr[\text{$\subspace{U}$ is realizable}] \geq 1 - 7.8\sqrt{n}e^{-n/24}.\]

\section{Tractable block diagonal and low-rank decompositions and related problems}
\label{sec:block}

In this section we generalize our results to the analogue of MTFA for
\emph{block}-diagonal and low-rank decompositions. Mimicking our earlier
development, we relate the analysis of this variant of MTFA to the facial
structure of a variant of the elliptope and a generalization of the ellipsoid
fitting problem. The key point is that these problems all possess additional
symmetries that, once taken into account, essentially allow us to reduce our
analysis to cases already considered in Sections~\ref{sec:links} and~\ref{sec:suff}.

Throughout this section, let $\mathcal{P}$ be a fixed partition of
$\{1,2,\ldots,n\}$.  We say a matrix is $\mathcal{P}$-block-diagonal if it is
zero except for the principal submatrices indexed by the elements of
$\mathcal{P}$. We denote by $\blkdiag_{\mathcal{P}}$ the map that takes an $n\times n$ 
matrix and maps it to the principal submatrices
indexed by $\mathcal{P}$. Its adjoint, denoted $\blkdiag_{\mathcal{P}}^*$,
takes a tuple of symmetric matrices $(X_\I)_{\I\in \mathcal{P}}$ and produces
an $n\times n$ matrix that is $\mathcal{P}$-block diagonal with blocks given by the $X_{\I}$. 

We now describe the analogues of MTFA, ellipsoid fitting, and
the problem of determining the facial structure of the elliptope.

\paragraph{Block minimum trace factor analysis}
If $X = B^\star + L^\star$ where $B^\star$ is $\mathcal{P}$-block-diagonal and
$L^\star \psd 0$ is low rank, the obvious analogue of MTFA is the semidefinite
program
\begin{equation}
    \label{eq:bmtfa}
    \minimize_{B,L} \; \tr(L) \quad\text{subject to}\quad
    \left\{ \begin{array}{rcl} X \!\!\!\! & = & \!\!\!\! B + L\\
            L \!\!\!\! & \psd \!\!\!\! & 0\\
            B\!\!\!\! & \text{is} &\!\!\!\! \text{$\mathcal{P}$-block-diagonal}
        \end{array}\right.
\end{equation}
which we call \emph{block minimum trace factor analysis} (BMTFA). 
\begin{definition}
    \label{def:rBMTFA}
    A subspace $\subspace{U}$ of $\R^n$ is \emph{recoverable by BMTFA} if for every $B^\star$ that 
    is $\mathcal{P}$-block-diagonal and every positive semidefinite $L^\star$ with column 
    space $\subspace{U}$, $(B^\star,L^\star)$ is the unique
    optimum of BMTFA with input $X = B^\star + L^\star$. 
\end{definition}

\paragraph{Faces of the $\mathcal{P}$-elliptope}
Just as MTFA is related to the facial structure of the elliptope, BMTFA is related to the facial
structure of the spectrahedron 
\[ \elliptope_{\mathcal{P}} = \{Y \psd 0: \blkdiag_{\mathcal{P}}(Y) = (I,I,\ldots,I)\}.\]
We refer to $\elliptope_\mathcal{P}$ as the $\mathcal{P}$-elliptope.
We extend the definition of a realizable subspace to this context.
\begin{definition}
    \label{def:Prealizable}
    A subspace $\subspace{U}$ of $\R^n$ is \emph{$\mathcal{P}$-realizable} if there is some $Y\in
    \elliptope_{\mathcal{P}}$ such that $\nullsp(Y) \supseteq \subspace{U}$.
\end{definition}

\paragraph{Generalized ellipsoid fitting}

To describe the $\mathcal{P}$-ellipsoid fitting problem 
we first introduce some convenient notation. If $\I\subset [n]$ we write
\begin{equation}
    \label{eq:SI}
    S^{\I} = \{x\in \R^n: \|x\|_2 = 1,\;\text{$x_j = 0$ if $j\notin\I$}\}
\end{equation}
for the intersection of the unit sphere with the coordinate subspace indexed by $\I$. 

Suppose $v_1,v_2,\ldots,v_n\in \R^k$ is a collection of points and $V$ is the $k\times n$ 
matrix with columns given by the $v_i$. Noting that $S^{\{i\}} = \{-e_i,e_i\}$, and thinking of $V$
as a linear map from $\R^n$ to $\R^k$, we see that 
the ellipsoid fitting problem is to find an ellipsoid in $\R^k$ with boundary containing 
$\cup_{i\in [n]} V(S^{\{i\}})$, i.e.~the collection of points $\pm v_1,\ldots,\pm v_n$. 
The $\mathcal{P}$-ellipsoid fitting problem is then to find an ellipsoid in 
$\R^k$ with boundary containing $\cup_{\I\in \mathcal{P}} V(S^{\I})$, i.e.~the collection of
ellipsoids $V(S^{\I})$. 

The generalization of the ellipsoid fitting property of a subspace is as follows.
\begin{definition}
    \label{def:Pefp}
    A subspace $\subspace{V}$ of $\R^n$ has the \emph{$\mathcal{P}$-ellipsoid fitting property}
    if there is a $k\times n$ matrix $V$ with row space $\subspace{V}$
    such that there is a centered ellipsoid in $\R^k$ with boundary containing
    $\cup_{\I\in \mathcal{P}} V(S^{\I})$.
\end{definition}

\subsection{Relating the generalized problems}
The facial structure of the $\mathcal{P}$-elliptope, BMTFA, and the
$\mathcal{P}$-ellipsoid fitting problem are related by the following
result, the proof of which is omitted as it is almost identical to that of Proposition~\ref{prop:links}.
\begin{proposition}
    \label{prop:block-links}
    Let $\subspace{U}$ be a subspace of $\R^n$. 
    Then the following are equivalent:
    \begin{remunerate}
    \item $\subspace{U}$ is recoverable by BMTFA.
        \item $\subspace{U}$ is $\mathcal{P}$-realizable.
        \item $\subspace{U}^\perp$ has the $\mathcal{P}$-ellipsoid fitting property.
    \end{remunerate}
\end{proposition}

The following lemma is the analogue of Lemma~\ref{lem:failure}. It describes certificates that
a subspace $\subspace{U}$ is not $\mathcal{P}$-realizable. Again the proof is almost identical to 
that of Lemma~\ref{lem:failure} so we omit it.
\begin{lemma}
    \label{lem:block-failure}
    A subspace $\subspace{U}$ of $\R^n$ is \emph{not} $\mathcal{P}$-realizable if and only if 
    there is a $\mathcal{P}$-block-diagonal matrix $B$ such that $\tr(B) > 0$ and 
    $v^TBv \leq 0$ for all $v\in \subspace{U}^\perp$. 
\end{lemma}

For the sake of brevity, in what follows we only discuss the problem of whether $\subspace{U}$ is $\mathcal{P}$-realizable
without explicitly translating the results into the context of the other two problems.

\subsection{Symmetries of the $\mathcal{P}$-elliptope}

We now consider the symmetries of the $\mathcal{P}$-elliptope. Our motivation for doing 
so is that it allows us to partition subspaces into classes for which either all 
elements are $\mathcal{P}$-realizable or none of the elements are $\mathcal{P}$-realizable.

It is clear that the $\mathcal{P}$-elliptope is invariant under conjugation 
by $\mathcal{P}$-block-diagonal orthogonal matrices. Let $G_{\mathcal{P}}$ denote
this subgroup of the group of $n\times n$ orthogonal matrices. There is a natural action of $G_{\mathcal{P}}$ on 
subspaces of $\R^n$ defined as follows. If $P\in G_{\mathcal{P}}$ and 
$\subspace{U}$ is a subspace of $\R^n$ then $P\cdot \subspace{U}$
is the image of the subspace $\subspace{U}$ under the map $P$.
(It is straightforward to check that this is a well defined group action.) If there exists 
some $P\in G_{\mathcal{P}}$ such that $P\cdot \subspace{U} = \subspace{U}'$ then 
we write $\subspace{U} \sim \subspace{U}'$ and say that $\subspace{U}$ and $\subspace{U}'$ are \emph{equivalent}.
We care about this equivalence relation on subspaces because the property of being $\mathcal{P}$-realizable 
is really a property of the corresponding equivalence classes.
\begin{proposition}
    \label{prop:real-equiv}
    Suppose $\subspace{U}$ and $\subspace{U}'$ are subspaces of $\R^n$. If $\subspace{U}\sim \subspace{U}'$ then
    $\subspace{U}$ is $\mathcal{P}$-realizable if and only if $\subspace{U}'$ is $\mathcal{P}$-realizable.
\end{proposition}

\begin{proof}
    If $\subspace{U}$ is $\mathcal{P}$-realizable there is $Y\in
    \elliptope_{\mathcal{P}}$ such that $Yu = 0$ for all $u\in \subspace{U}$.
    Suppose $\mathcal{U}' = P\cdot \mathcal{U}$ for some $P\in G_{\mathcal{P}}$
    and let $Y' = PYP^T$. Then $Y'\in \elliptope_{\mathcal{P}}$ and $Y'(Pu) =
    (PYP^T)(Pu) = 0$ for all $u\in \subspace{U}$.  By the definition of
    $\subspace{U}'$ it is then the case that $Y'u' = 0$ for all $u'\in
    \subspace{U}'$. Hence $\subspace{U}'$ is $\mathcal{P}$-realizable. The
    converse clearly also holds.
\end{proof}

\subsection{Exploiting symmetries: relating realizability and $\mathcal{P}$-realizability}
\label{sec:sym-real}
For a subspace of $\R^n$, we now consider how the notions of $\mathcal{P}$-realizability and realizability (i.e.~$[n]$-realizability)
relate to each other. Since $\elliptope_{\mathcal{P}} \subset \elliptope_n$, if $\subspace{U}$ is 
$\mathcal{P}$-realizable, it is certainly also realizable. While the converse does not hold, we can establish the
following partial converse, which we subsequently use to extend our analysis from Sections~\ref{sec:links} and~\ref{sec:suff}
to the present setting. 
\begin{theorem}
    \label{thm:sym-real}
    A subspace $\subspace{U}$ of $\R^n$ is $\mathcal{P}$-realizable if and only if
    $\subspace{U}'$ is realizable for every $\subspace{U}'$ such that $\subspace{U}' \sim \subspace{U}$.
\end{theorem} 

\begin{proof}
    We note that one direction of the proof is obvious since $\mathcal{P}$-realizability implies realizability.
It remains to show that if $\subspace{U}$ is not $\mathcal{P}$-realizable
then there is some $\subspace{U}'$ equivalent to $\subspace{U}$ that is not realizable. 

Recall from Lemma~\ref{lem:block-failure} that if $\subspace{U}$ is not $\mathcal{P}$-realizable
there is some $\mathcal{P}$-block-diagonal $X$ with positive trace such that $v^T X v \leq 0$ 
for all $v\in \subspace{U}^\perp$. Since $X$ is $\mathcal{P}$-block-diagonal there is some 
$P\in G_{\mathcal{P}}$ such that $PXP^T$ is diagonal. Since conjugation by orthogonal matrices 
preserves eigenvalues, $\tr(PXP^T) = \tr(X) > 0$. Furthermore
$v^T(PXP^T)v = (P^Tv)^TX(P^Tv) \leq 0$ for all $P^Tv \in \subspace{U}^\perp$. Hence 
$w^T(PXP^T)w \geq 0$ for all $w\in P\cdot\subspace{U}^\perp = (P\cdot \subspace{U})^{\perp}$. 
By Lemma~\ref{lem:failure}, $PXP^T$ is a certificate that $P\cdot \subspace{U}$ is not realizable, completing the proof.
\end{proof}

The power of Theorem~\ref{thm:sym-real} lies in its ability to turn any condition for a subspace to be realizable 
into a condition for the subspace to be $\mathcal{P}$-realizable by appropriately symmetrizing 
the condition with respect to the action of $G_{\mathcal{P}}$. We now illustrate this approach by 
generalizing Theorem~\ref{thm:DP} and our coherence based condition (Theorem~\ref{thm:main})
for a subspace to be $\mathcal{P}$-realizable. In each case we first define an appropriately symmetrized 
version of the original condition. The natural symmetrized version of the notion of balance is as follows.
\begin{definition}
    \label{def:Pbal}
    A vector $u\in \R^n$ is \emph{$\mathcal{P}$-balanced} if for all $\I\in \mathcal{P}$
    \[ \|u_{\I}\|_2 \leq \sum_{\J\in \mathcal{P}\setminus\{\I\}} \|u_{\J}\|_2.\]
\end{definition}

We next define the appropriately symmetrized analogue of coherence. Just as coherence measures 
how far a subspace is from any one-dimensional coordinate subspace, 
$\mathcal{P}$-coherence measures how far a subspace is from any of the coordinate subspaces indexed by 
elements of $\mathcal{P}$. 
\begin{definition}
    \label{def:Pcoh}
    The $\mathcal{P}$-coherence of a subspace $\subspace{U}$ of $\R^n$ is 
\[ \mu_{\mathcal{P}}(\subspace{U}) = \max_{\I\in \mathcal{P}}\max_{x\in S^{\I}} \|P_{\subspace{U}}x\|_2^2.\]
\end{definition}

Just as the coherence of $\subspace{U}$ can be computed by taking the maximum diagonal element of 
$P_{\subspace{U}}$, it is straightforward to veify that the $\mathcal{P}$-coherence of $\subspace{U}$ 
can be computed by taking the maximum of the spectral norms of the principal submatrices 
$[P_{\subspace{U}}]_{\I}$ indexed by $\I\in \mathcal{P}$.

We now use Theorem~\ref{thm:sym-real} to establish the natural generalization of Theorem~\ref{thm:DP}.
\begin{corollary}
    \label{cor:Pbal-Preal}
    If a subspace $\subspace{U}$ of $\R^n$ is $\mathcal{P}$-realizable then every element of $\subspace{U}$ is 
    $\mathcal{P}$-balanced. If $\subspace{U} = \spn\{u\}$ is one dimensional then 
    $\subspace{U}$ is $\mathcal{P}$-realizable if and only if $u$ is $\mathcal{P}$-balanced.
\end{corollary}
\begin{proof}    
    If there is $u\in \subspace{U}$ that is not $\mathcal{P}$-balanced then there is 
    $P\in G_{\mathcal{P}}$ such that $Pu$ is not balanced (choose $P$ so that it rotates 
    each $u_{\I}$ until it has only one non-zero entry). But then $P\cdot \subspace{U}$ is 
    not realizable and so $\subspace{U}$ is not $\mathcal{P}$-realizable.

    For the converse, we first show that if a vector is $\mathcal{P}$-balanced then it is balanced. 
    Let $\I\in \mathcal{P}$, and consider $i\in \I$. Then since $u$ is $\mathcal{P}$-balanced,
    \[ 2|u_i| \leq 2 \|u_{\I}\|_2 \leq \sum_{\J\in \mathcal{P}} \|u_{\J}\|_2 \leq \sum_{i=1}^{n}|u_i|\]
    and so $u$ is balanced.

    Now suppose $\subspace{U} = \spn\{u\}$ is one dimensional and $u$ is $\mathcal{P}$-balanced.
    Since $u$ is $\mathcal{P}$-balanced it follows that $Pu$ is $\mathcal{P}$-balanced (and hence balanced)
    every $P\in G_{\mathcal{P}}$. Then by Theorem~\ref{thm:DP}
    $\spn\{Pu\}$ is realizable for every $P\in G_{\mathcal{P}}$. Hence by Theorem~\ref{thm:sym-real}, 
    $\subspace{U}$ is $\mathcal{P}$-realizable. 
\end{proof}

Similarly, with the aid of Theorem~\ref{thm:sym-real} we can write down a $\mathcal{P}$-coherence-threshold
condition that is a sufficient condition for a subspace to be $\mathcal{P}$-realizable. 
The following is a natural generalization of Theorem~\ref{thm:main}.
\begin{corollary}
    \label{cor:Pcoh-Preal}
    If $\mu_{\mathcal{P}}(\subspace{U}) < 1/2$ then $\subspace{U}$ is $\mathcal{P}$-realizable.
\end{corollary}
\begin{proof}
    By examining the constraints in the variational definitions of $\mu(\subspace{U})$ and 
    $\mu_{\mathcal{P}}(\subspace{U})$ we see that $\mu(\subspace{U}) \leq \mu_{\mathcal{P}}(\subspace{U})$.
    Consequently if $\mu_{\mathcal{P}}(\subspace{U}) < 1/2$ it follows from Theorem~\ref{thm:main}
    that $\subspace{U}$ is realizable. 
    Since $\mu_{\mathcal{P}}$ is invariant under the action of $G_{\mathcal{P}}$ on subspaces we
    can apply Theorem~\ref{thm:sym-real} to complete the proof.
\end{proof}

\section{Conclusions}
We established a link between three problems of independent interest: 
deciding whether there is a centered ellipsoid passing through a
collection of points, understanding the structure of the faces of the
elliptope, and deciding which pairs of diagonal and low rank-matrices can be
recovered from their sum using a tractable semidefinite-programming-based
heuristic, namely minimum trace factor analysis.  We provided a simple sufficient
condition, based on the notion of the coherence of a subspace, which ensures the
success of minimum trace factor analysis, and showed that this is the best
possible coherence-threshold-type sufficient condition for this problem.  We
provided natural generalizations of our results to the problem of analyzing
tractable block-diagonal and low-rank decompositions, showing how the
symmetries of this problem allow us to reduce much of the analysis to the
original diagonal and low-rank case.

Our results suggest both the power and the limitations of using `coarse'
properties of a subspace such as coherence to gain understanding of the faces of the elliptope
(and related problems). The power of results based on such properties is that they do not have explicit
dimension-dependence, unlike previous results on the faces of the elliptope.  At
the same time, the lack of explicit dimension dependence typically yields conservative sufficient
conditions for high-dimensional problems. It would be interesting to find a
hierarchy of coherence-like conditions that provide less conservative
sufficient conditions for higher dimensional problem instances.

\appendix

\section{Additional proofs}
\label{app:pfs}

\subsection{Proof of Lemma~\ref{lem:bal-ef}}
We first establish Lemma~\ref{lem:bal-ef} which gives an interpretation of the balance condition in terms of 
ellipsoid fitting.
\begin{proof}
    The proof is a fairly straightforward application of linear programming duality. Throughout let $V$ be the 
    $k\times n$ matrix with columns given by the $v_i$.  
    The point $v_i\in \R^k$ is on the boundary of the convex hull of $\pm v_1,\ldots,\pm v_n$
    if and only if there exists $x\in \R^k$ such that $\langle x,v_i\rangle = 1$ and 
    $|\langle x,v_j\rangle| \leq 1$ for all $j\neq i$. Equivalently, the following linear program (which depends on $i$)
    is feasible
    \begin{equation}
            \label{eq:lp-primal}
        \minimize_x \; \langle 0,x\rangle \quad\text{subject to}\quad 
        \left\{\begin{array}{rcl} v_i^Tx \!\!\!\! & = & \!\!\!\! 1\\
                |v_j^Tx| \!\!\!\! & \leq & \!\!\!\! 1\;\;\text{for all $j\neq i$.}\end{array}\right.
        \end{equation}
    Suppose there is some $i$ such that $v_i$ is in the interior of $\text{conv}\{\pm v_1,\ldots,\pm v_n\}$.
    Then \eqref{eq:lp-primal} is not feasible so the dual linear program (which depends on $i$) 
    \begin{equation}
        \maximize_u \; u_i - \sum_{j\neq i} |u_j| \quad\text{subject to}\quad
        \begin{array}{rcl} V u \!\!\!\! & = & \!\!\!\! 0\end{array}\label{eq:lp-dual}
           \end{equation}  
    is unbounded. This is the case if and only if there is some $u$ in the nullspace of $V$
    such that $u_i > \sum_{j\neq i} |u_j|$. If such a $u$ exists, then it is certainly the 
    case that $|u_i| \geq u_i > \sum_{j\neq i} |u_j|$ and so $u$ is not balanced. 
    
    Conversely if $u$ is in the nullspace of $V$ and $u$ is not balanced then either $u$ or $-u$
    satisfies $u_i > \sum_{j\neq i} |u_j|$ for some $i$. Hence the linear program 
    \eqref{eq:lp-dual} associated with the index $i$ is unbounded and so the corresponding  
    linear program \eqref{eq:lp-primal} is infeasible. It follows that $v_i$ is in the interior of the convex
    hull of $\pm v_1,\ldots,\pm v_n$. 
\end{proof}
\subsection{Completing the proof of Theorem~\ref{thm:main}}
We now complete the proof of Theorem~\ref{thm:main} by establishing the following result about the 
existence of a non-negative solution to the linear system~\eqref{eq:lin-sys}.
\begin{lemma}
    \label{lem:linsys}
    If $\mu(\subspace{U}) < 1/2$ then there is $\lambda \geq 0$ such that 
    \begin{equation}
        \label{eq:linsys}
        \diag\left(P_{\subspace{U}^\perp}\diag^*(\lambda)P_{\subspace{U}^\perp}\right) = \ones.
    \end{equation}
\end{lemma}

\begin{proof}
     We note that the linear system \eqref{eq:linsys} can be written as 
     $P_{\subspace{U}^\perp}\circ P_{\subspace{U}^\perp}\lambda = \ones$
    where $\circ$ denotes the entry-wise product of matrices.
    As such, we need to show that $P_{\subspace{U}^\perp}\circ P_{\subspace{U}^\perp}$ is invertible and
    $(P_{\subspace{U}^\perp}\circ P_{\subspace{U}^\perp})^{-1}\ones \geq 0$. To do so, we appeal to the 
    following (slight restatement) of a theorem of Walters~\cite{walters1969nonnegative} regarding positive solutions
    to certain linear systems.
\begin{theorem}[Walters \cite{walters1969nonnegative}]
    \label{thm:walters}
	Suppose $A$ is a square matrix with non-negative entries and positive diagonal entries. 
    Let $D$ be a diagonal matrix with $D_{ii} = A_{ii}$ for all $i$. If $y>0$ and $2y - AD^{-1}y > 0$ then 
    $A$ is invertible and $A^{-1}y > 0$.
\end{theorem}

In our case we take $A = P_{\subspace{U}^\perp}\circ P_{\subspace{U}^\perp}$ and $y = \ones$ in Theorem~\ref{thm:walters}. It is
clear that $P_{\subspace{U}^\perp} \circ P_{\subspace{U}^\perp}$ is entry-wise non-negative. 
Furthermore 
$[P_{\subspace{U}^\perp}]_{ii} = 1-[P_\subspace{U}]_{ii} > 1-\mu(\subspace{U}) > 1/2$
and so $D_{ii} = [P_{\subspace{U}^\perp}\circ P_{\subspace{U}^\perp}]_{ii} > 1/4$.
It then remains to show that $P_{\subspace{U}^\perp}\circ~P_{\subspace{U}^\perp}~D^{-1}\ones<2\ones$. Consider the $i$th
such inequality, and observe that
\begin{align*}
    [P_{\subspace{U}^\perp}\circ P_{\subspace{U}^\perp} D^{-1}\ones]_i 
    & = \left(P_{\subspace{U}^\perp}D^{-1}P_{\subspace{U}^\perp}\right)_{ii}\\
    & = \left(P_{\subspace{U}^\perp}D_{ii}^{-1}e_ie_i^TP_{\subspace{U}^\perp}\right)_{ii} + 
    \left(P_{\subspace{U}^\perp}(D^{-1} - D_{ii}^{-1}e_ie_i^T)P_{\subspace{U}^\perp}\right)_{ii}\\
    & \leq 1 + \max_{j\in[n]}D_{jj}^{-1} \left(P_{\subspace{U}^\perp}(I-e_ie_i^T)P_{\subspace{U}^\perp}\right)_{ii}\\
    & < 1 + 4[P_{\subspace{U}^\perp}]_{ii} - 4[P_{\subspace{U}^\perp}]_{ii}^2\\
    & = 2 - 4([P_{\subspace{U}^\perp}]_{ii} - 1/2)^{2}\\
    & \leq 2
\end{align*}
where we have used the assumption that $[P_{\subspace{U}^\perp}]_{ii}>1/2$ for all $i$
and the fact that $P_{\subspace{U}^\perp}^{2} = P_{\subspace{U}^\perp}$. Applying Walters's theorem 
completes the proof.
\end{proof}

\subsection{Proof of Proposition~\ref{prop:rand}}
We now establish Proposition~\ref{prop:rand}, giving a bound on the probability that a suitably 
random subspace is realizable by bounding the probability that it has coherence strictly bounded above by $1/2$.
\begin{proof}
    It suffices to show that $\|P_{\subspace{U}} e_i\|^2 \leq 
    (1-2\epsilon)(1/2-\epsilon) = 1/2 - 2\epsilon^2 < 1/2$ for all $i$ with high probability. 
        The main observation we use is that if $\subspace{U}$ is a random
        $r$ dimensional subspace of $\R^n$ and $x$ is any fixed vector with
        $\|x\|=1$ then $\|P_{\subspace{U}} x\|^2 \sim \beta(r/2,(n-r)/2)$ where $\beta(p,q)$ denotes the
        beta distribution \cite{frankl1990some}. In the case where $r = (1/2 - \epsilon)n$, using a tail bound 
        for $\beta$ random variables \cite{frankl1990some} we see that if $x\in \R^n$ is fixed 
        and $r > 3/\epsilon^2$ then 
        \[ \Pr[ \|P_\subspace{U} x\|^2 \geq (1+2\epsilon)(1/2 -\epsilon)] < 
            \frac{1}{a_\epsilon}\frac{1}{(\pi(1/4 - \epsilon^2))^{1/2}}
            n^{-1/2}e^{-a_\epsilon k}\]
        where $a_\epsilon = \epsilon - 4\epsilon^2/3$.
        Taking a union bound over $n$ events, as long as $r > 3/\epsilon^2$ 
        \begin{align*}
            \Pr\left[\mu(\subspace{U}) \geq 1/2\right] & \leq \Pr\left[\|P_{\subspace{U}} e_i\|^2 \geq (1-2\epsilon)(1/2-\epsilon)
                \;\; \text{for some $i\in [n]$}\right]\\
            & \leq n\cdot \frac{1}{a_\epsilon(\pi(1/4 - \epsilon^2))^{1/2}} n^{-1/2} e^{- a_\epsilon k} =
            \bar{c}n^{1/2}e^{-\tilde{c}n}
        \end{align*}
        for appropriate positive constants $\bar{c}$ and $\tilde{c}$.
\end{proof}

\subsection*{Acknowledgements} The authors would like to thank Prof.~Sanjoy Mitter for helpful discussions.

\bibliographystyle{siam}
\bibliography{dlr-bib}

\end{document}